\documentclass[]{article}
\usepackage{amssymb,amsfonts,amsmath}
\usepackage{latexsym}
\usepackage{xcolor}

\newtheorem{thm}{Theorem}
\newtheorem{pro}[thm]{Proposition}
\newtheorem{lem}[thm]{Lemma}

\newcommand{\prf}{\noindent{\it Proof.} }
\newcommand{\cbdo}{\hfill$\Box$}
\newcommand{\ind}{\mathbf 1}

\begin{document}

\title{The Erd\H{o}s four-edge intersection problem}
\author{Andrzej \.{Z}ak\thanks{The author was partially supported
by the Polish Ministry of Science and Higher Education.}\\
\small{AGH University of Krakow, Poland}}
\maketitle

\begin{abstract}
For an $n$-vertex graph $G$ and a permutation $\sigma$ of its vertex set,
let $\sigma(G)$ denote the corresponding relabelling of $G$, and put
\[
 I_G(\sigma)=|E(G)\cap E(\sigma(G))|.
\]
Let $f(n,k)$ be the minimum number of edges in an $n$-vertex graph for
which $I_G(\sigma)\geq k$ for every $\sigma$.  In 1977 Erd\H{o}s asked
whether $f(n,4)=2n-4$, observing that $K_{2,n-2}$ gives the upper bound.
We prove that, for all sufficiently large $n$,
\[
 f(n,4)=2n-4.
\]
Equivalently, every sufficiently large $n$-vertex graph with at most
$2n-5$ edges has a relabelling with at most three common edges.

Our proof is inspired by the recent work of Fang and Hou on the
Erd\H{o}s--Mullin five-edge intersection problem and builds on their
core--buffer and absorption framework.  The main additional ingredients
are a growing high-degree core $C$ satisfying
\[
 |C|\Delta(G-C)=o(n),
\]
and a rigidity analysis of the equality case in the relevant first-moment
estimate.  This analysis shows that the only core--buffer configuration
forcing four local common edges is of $K_{2,|C|}$ type; the strict bound
$e(G)\leq2n-5$ then supplies a defect which breaks this configuration.
\end{abstract}

\section{Introduction}

Throughout the paper graphs are finite and simple.  If $\sigma$ is a
permutation of $V(G)$, then
\[
 E(\sigma(G))=\{\sigma(u)\sigma(v):uv\in E(G)\},
\]
and
\[
 I_G(\sigma)=|E(G)\cap E(\sigma(G))|.
\]
For $k\geq1$, let
\[
 f(n,k):=\min\bigl\{|E(G)|:|V(G)|=n\text{ and }
 I_G(\sigma)\geq k\text{ for every permutation }\sigma\bigr\}.
\]
Thus every $n$-vertex graph with fewer than $f(n,k)$ edges has a
relabelling sharing at most $k-1$ edges with the original graph.

Erd\H{o}s~\cite{Erdos} considered the small values of this function.
After discussing the already settled cases $k\leq3$, he left the cases
$k=4$ and $k=5$ open.  For $k=4$ he observed that
\[
 f(n,4)\leq2n-4
\]
by taking $K_{2,n-2}$, and asked whether equality holds.  Immediately
afterwards he recorded Mullin's construction giving $f(n,5)\leq2n-2$
and suggested equality there as well.  The latter problem has recently
been solved by Fang and Hou~\cite{FH}, who proved
$f(n,5)=2n-2$ for all sufficiently large $n$.

The minimum-intersection viewpoint was placed in a broader framework by
Mullin, Roy, and Schellenberg~\cite{MRS}.  We use the language of
near-packings from~\cite{zak1}.  For $r\geq0$, let $\mathcal E_r$ be the
family of graphs with at most $r$ edges, and let $\mu(n,\mathcal E_r)$ be
the largest integer $m$ such that every $n$-vertex graph with at most
$m$ edges has an $\mathcal E_r$-near-packing.  Directly from the
definitions,
\begin{equation}\label{relation-f-mu}
 f(n,k)=\mu(n,\mathcal E_{k-1})+1.
\end{equation}
Hence the four-edge problem is the determination of
$\mu(n,\mathcal E_3)$.

Our main result is the following.

\begin{thm}\label{main}
For all sufficiently large $n$,
\[
 f(n,4)=2n-4.
\]
Equivalently, every sufficiently large $n$-vertex graph $G$ with
\[
 e(G)\leq2n-5
\]
has a permutation $\sigma$ such that $I_G(\sigma)\leq3$.
\end{thm}

The upper bound is Erd\H{o}s's graph $K_{2,n-2}$.  Indeed, if the
two-vertex parts of two copies are disjoint, the four edges between
these parts are common.  If the two parts meet in one vertex, there are
$n-2$ common edges, and if they coincide there are $2n-4$ common edges.

Our proof builds on the core--buffer and absorption framework developed
by Fang and Hou~\cite{FH} for the five-edge problem.  In particular, we
use their sparse-permutation list-packing lemma, their degree-deficit
mechanism for producing low-degree buffer vertices, and their absorption
scheme for the external neighbours of the buffer.  The main new
difficulty is that the first-moment calculation used in~\cite{FH} gives
four common edges exactly at the critical density $e(C,B)=2|C|$, whereas
here one has to save one further edge.  We exploit the strict inequality
\[
 e(C,B)<2|C|
\]
to obtain expectation strictly smaller than four and then analyse the
equality case.  The resulting rigidity lemma forces a
$K_{2,|C|}$-type configuration; the strict global bound
$e(G)\leq2n-5$ then supplies a local defect which breaks this rigidity.
A second new feature is that we replace the fixed finite core of
Fang--Hou by a growing core satisfying
\[
 |C|=o(n),\qquad |C|\Delta(G-C)=o(n),
\]
which is obtained directly from the degree ordering and allows the
absorption argument to be run without a diagonal limit of the normalized
degree sequence.  In the almost-universal case we likewise refine the
corresponding Fang--Hou construction so as to eliminate one additional
possible common edge.

We make the dependence on~\cite{FH} explicit throughout.  Whenever an
auxiliary statement can be used in exactly the form proved there, we
quote it and do not repeat its proof.  For the sake of completeness, we
retain proofs of the auxiliary variants needed in our argument even when
their constructions are close to those in~\cite{FH}; this also makes
clear where the growing-core hypothesis or the passage from four common
edges to three requires an additional argument.  We also organize the three
extension steps through a single completion principle: after a local partial
bijection has been prescribed, all possible new common edges are encoded by
a forbidden-pair graph and the remainder is handled by list packing.

\section{Packing tools}

We use the following classical packing theorem.

\begin{thm}[\cite{BE,BS,SS}]\label{pak1}
Every graph of order $N$ and size at most $N-2$ is packable.
\end{thm}

The same first-moment argument as in~\cite[Theorem~17]{zak1} gives the
following general estimate, which will also be useful below.

\begin{pro}\label{firstmoment}
For integers $N\geq3$ and $0\leq k<\binom N2$,
\[
 \mu(N,\mathcal E_k)\geq
 \left\lceil\sqrt{(k+1)\binom N2}\right\rceil-1.
\]
\end{pro}

\prf
Let $G$ be an $N$-vertex graph with $m$ edges, and choose a permutation
$\sigma$ uniformly at random.  Each edge of $G$ is mapped uniformly to
one of the $\binom N2$ pairs of vertices, and hence
\[
 \mathbb E I_G(\sigma)=\frac{m^2}{\binom N2}.
\]
If
\[
 m\leq\left\lceil\sqrt{(k+1)\binom N2}\right\rceil-1,
\]
then $m^2<(k+1)\binom N2$, so $\mathbb E I_G(\sigma)<k+1$.
Therefore some permutation satisfies $I_G(\sigma)\leq k$. \cbdo

In particular, taking $k=1$, every $N$-vertex graph with at most $N-1$
edges has a relabelling with at most one common edge.

We next recall list packing.  Let $F_1$ and $F_2$ be graphs on two
disjoint $N$-vertex sets, and let $Y$ be a bipartite graph between these
sets.  A bijection $\phi:V(F_1)\to V(F_2)$ is a \emph{list packing} if
it maps every edge of $F_1$ to a non-edge of $F_2$ and
$x\phi(x)\notin E(Y)$ for every $x\in V(F_1)$.  The edges of $Y$ are
forbidden vertex--image pairs.

\begin{thm}[Gy\H{o}ri--Kostochka--McConvey--Yager~\cite{GKMY}]
\label{listSS}
If
\[
 \Delta(F_1)\Delta(F_2)+\Delta(Y)<\frac N2,
\]
then $F_1,F_2$ have a list packing avoiding $Y$.
\end{thm}

The final completion is used when both graphs have linearly many
edges and all three maximum degrees are a small linear fraction of the
order.  We use the sparse permutation list-packing lemma of Fang and
Hou~\cite[Lemma~2.3]{FH} in exactly the form below, and therefore omit
its proof.

\begin{lem}[Fang--Hou, Lemma 2.3]\label{sparse}
There is an absolute constant $\varepsilon_0>0$ such that the following
holds for all sufficiently large $N$.  Let $F_1,F_2$ be $N$-vertex
graphs and let $Y$ be a bipartite graph between their vertex sets.  If
\[
 e(F_1),e(F_2)\leq3N
\]
and
\[
 \max\{\Delta(F_1),\Delta(F_2),\Delta(Y)\}
 \leq\varepsilon_0N,
\]
then $F_1,F_2$ have a list packing avoiding $Y$.
\end{lem}

All subsequent extension arguments use the following elementary completion
principle.  For a partial bijection $p:A\to B$ between the source and target
copies of a graph $F$, put
\[
 U=V(F)\setminus A,\qquad W=V(F)\setminus B
\]
and let
\[
 I_F(p)=\bigl|\{aa'\in E(F[A]):p(a)p(a')\in E(F)\}\bigr|.
\]
Define a bipartite graph $Y_p$ between the source copy of $U$ and the target
copy of $W$ by declaring $xy\in E(Y_p)$ whenever there is $a\in A$ such that
\[
 ax\in E(F)\qquad\text{and}\qquad p(a)y\in E(F).
\]
Thus $Y_p$ records precisely the vertex--image pairs which would create a
new common edge together with the already prescribed part.

\begin{lem}\label{completion}
If $F[U],F[W]$ have a list packing $\phi:U\to W$ avoiding $Y_p$, then
$p\cup\phi$ is a permutation of $V(F)$ and
\[
 I_F(p\cup\phi)=I_F(p).
\]
\end{lem}

\prf
Edges with both ends in $A$ contribute exactly $I_F(p)$.  Edges with both
ends in $U$ create no common edge because $\phi$ is a packing of
$F[U]$ with $F[W]$.  Finally, if $ax\in E(F)$ with $a\in A$ and $x\in U$,
then $x\phi(x)\notin E(Y_p)$ implies $p(a)\phi(x)\notin E(F)$.  Hence no
new common edge is created. \cbdo

A bounded prescribed part will be needed in the almost-universal case.  The
next extension lemma is extracted from the list-packing construction in the
proof of Fang and Hou~\cite[Lemma~3.1]{FH}.  We state it in a form which
allows a bounded partial near-packing and preserves all common edges already
present in the prescribed part.

\begin{lem}\label{prescribed}
Let $(F_N)$ be a sequence of graphs with $|V(F_N)|=N$, and put
\[
 H_N=\{v\in V(F_N):d_{F_N}(v)>2\}.
\]
Suppose, as $N\to\infty$, that
\[
 \Delta(F_N)=o(N),\qquad |H_N|=o(N),\qquad
 \sum_{v\in H_N}d_{F_N}(v)=o(N).
\]
Then, for every fixed $r$, every partial bijection $p:A\to B$ with
$|A|=|B|\leq r$ extends, for all sufficiently large $N$, to a permutation
$\sigma$ of $V(F_N)$ satisfying
\[
 I_{F_N}(\sigma)=I_{F_N}(p).
\]
\end{lem}

\prf
Write $F=F_N$ and $H=H_N$.  Enlarging $H$ by the bounded set
$A\cup B$ does not change the assumptions: the additional contribution
to the degree sum is $O(\Delta(F))=o(N)$.  Thus we may assume
$A\cup B\subseteq H$.

We first extend $p$ to all source vertices of $H$.  Suppose that
$h\in H$ has not yet been assigned an image.  We choose its image
$y$ outside $H$, so that $d_F(y)\leq2$, and require that no new common
edge be created.

There are only $o(N)$ forbidden choices for $y$.  Indeed, the initially
prescribed pairs forbid at most
\[
 r\Delta(F)
\]
vertices.  Moreover, if $h'$ is a previously assigned neighbour of
$h$, then the image of $h'$ lies outside $H$ and hence has degree at
most two; consequently all previously assigned neighbours of $h$
forbid at most
\[
 2d_F(h)
\]
further images.  Finally, at most $|H|+O(1)$ images have already been
used.  Thus the total number of excluded candidates is at most
\[
 r\Delta(F)+2d_F(h)+|H|+O(1)=o(N).
\]
Since $V(F)\setminus H$ has $N-o(N)$ vertices, the required image can
be chosen greedily.  Continuing in this way, we assign images to all
source vertices of $H$ without creating any new common edge.

Some vertices of $H$ may still be unused on the target side.  Let
$t$ be their number; clearly $t\leq |H|=o(N)$.  Choose an independent
set $R$ of size $t$ outside
\[
 H\cup N_F(H)
\]
and outside the vertices already used in the construction.  Such a set
exists: by assumption
\[
 |H|+\sum_{v\in H}d_F(v)=o(N),
\]
so the available set has $N-o(N)$ vertices, and its induced graph has
maximum degree at most two.  Hence a greedy choice of $t=o(N)$
independent vertices is possible.

Map $R$ bijectively onto the unused target vertices of $H$.  This
creates no new common edge among the prescribed vertices.  Indeed,
every previously prescribed source vertex belongs to $H$, while
$R$ is anticomplete to $H$, and the independence of $R$ prevents a
conflict between two newly prescribed pairs.  Possible conflicts
between a vertex of $R$ and a vertex whose image will be chosen later
are not required to be excluded at this stage; they are precisely
recorded by the forbidden graph in the completion step below.

Let $p'$ be the resulting partial bijection.  We have therefore
\[
 I_F(p')=I_F(p).
\]

Let $U$ and $W$ be the remaining source and target vertices,
respectively, and let $Y_{p'}$ be the forbidden graph defined in
Lemma~\ref{completion}.  Since all vertices of $H$ have been removed
from both the source and target sides,
\[
 \Delta(F[U]),\Delta(F[W])\leq2.
\]
Furthermore, a remaining source vertex has at most two prescribed
neighbours.  Each of their images has degree at most $\Delta(F)$, so
its degree in $Y_{p'}$ is at most $2\Delta(F)$.  Symmetrically, every
remaining target vertex also has degree at most $2\Delta(F)$ in
$Y_{p'}$.  Hence
\[
 \Delta(Y_{p'})\leq2\Delta(F)=o(N).
\]

If
\[
 N_0=|U|=|W|=N-o(N),
\]
then, for all sufficiently large $N$,
\[
 \Delta(F[U])\Delta(F[W])+\Delta(Y_{p'})
 \leq 4+2\Delta(F)<\frac{N_0}{2}.
\]
Theorem~\ref{listSS} therefore gives a list packing of $F[U]$ with
$F[W]$ avoiding $Y_{p'}$.  By Lemma~\ref{completion} it extends $p'$
to a permutation $\sigma$ of $V(F)$ without creating any further
common edge.  Consequently
\[
 I_F(\sigma)=I_F(p')=I_F(p),
\]
as required. \cbdo

\section{The core and the buffer}

The next elementary lemma replaces a diagonal limit of the normalized
degree sequence.

\begin{lem}\label{corechoice}
Let $G=G_n$ be an $n$-vertex graph with $e(G)\leq2n$, and let
$m=m(n)$ satisfy
\[
 m\longrightarrow\infty,
 \qquad
 m=o(n).
\]
Order the vertices so that
$d_1\geq d_2\geq\cdots\geq d_n$.  There is an integer
$1\leq k\leq m$ such that, if
\[
 C=\{v_i:i\leq k\text{ and }d_i\geq4\},
 \qquad c=|C|,
 \qquad D=\Delta(G-C),
\]
then
\[
 c=o(n),\qquad D=o(n),\qquad cD=o(n).
\]
Moreover, if $d_1\geq4$, then a vertex of maximum degree belongs to
$C$.
\end{lem}

\prf
Put
\[
 h_m=\sum_{j=1}^{m}\frac1j.
\]
Since $\sum_i d_i\leq4n$, there is $1\leq k\leq m$ such that
\[
 k d_{k+1}\leq\frac{4n}{h_m};
\]
otherwise
\[
 \sum_{j=1}^{m}d_{j+1}>
 \frac{4n}{h_m}\sum_{j=1}^{m}\frac1j=4n.
\]
Every vertex among the first $k$ which is not placed in $C$ has degree
at most three, and therefore
\[
 D\leq\max\{d_{k+1},3\}.
\]
Consequently
\[
 cD\leq k\max\{d_{k+1},3\}
 \leq\frac{4n}{h_m}+3m=o(n).
\]
This also gives $D=o(n)$, while $c\leq m=o(n)$.  The last assertion is
immediate. \cbdo

We now isolate a linear reservoir of low-degree vertices.  The following
degree-deficit estimate is a slight variant of Fang and
Hou~\cite[Lemma~4.1]{FH}.  For the sake of completeness we include the
short proof, since the present form uses the sharper bound $2n-5$ and a
core which is allowed to grow with $n$.

\begin{lem}\label{reservoir}
Let $G$ be an $n$-vertex graph with $e(G)\leq2n-5$, and let
$C\subseteq V(G)$ consist only of vertices of degree at least four.
Put
\[
 S_C=\sum_{v\in C}d(v)
\]
and
\[
 \mathcal B=\{v\notin C:d(v)\leq3\text{ and }d(v,C)\leq2\}.
\]
Then
\[
 |\mathcal B|\geq\frac{S_C}{6}-|C|+2.
\]
\end{lem}

\prf
Since $e(G)\leq2n-5$,
\[
 \sum_{v\in V(G)}(4-d(v))\geq10.
\]
Thus
\[
 \sum_{d(v)\leq3}(4-d(v))
 \geq
 \sum_{d(v)\geq5}(d(v)-4)+10
 \geq S_C-4|C|+10.
\]
Every vertex of degree at most three outside $\mathcal B$ has degree
three and all its neighbours in $C$.  If there are $r$ such vertices,
then
\[
 3r\leq S_C.
\]
Since each vertex of $\mathcal B$ contributes at most four to the
left-hand side above,
\[
 4|\mathcal B|+r\geq S_C-4|C|+10.
\]
Using $r\leq S_C/3$ gives the stated, slightly weaker, bound. \cbdo

The next lemma is the local reason why the number four is special.  Its
first-moment calculation is the one used by Fang and Hou in
Section~4.1 of~\cite{FH}; the analysis of the equality case is new and
is the point at which the four-edge problem departs from their argument.

\begin{lem}\label{corebuffer}
Let $C$ and $B$ be disjoint sets of the same size $s$, and suppose that
$B$ is independent and
\[
 d(b,C)\leq2\qquad(b\in B).
\]
Then one of the following holds.
\begin{enumerate}
\item There are bijections $\alpha:C\to B$ and $\beta:B\to C$ for which
at most three edges between $C$ and $B$ are mapped to edges between $B$
and $C$.
\item There is a two-element set $P\subseteq C$ such that
\[
 N(b)\cap C=P\qquad\text{for every }b\in B.
\]
\end{enumerate}
\end{lem}

\prf
Put $q=e(C,B)$.  Choose $\alpha$ and $\beta$ independently and uniformly
at random.  A fixed edge $cb\in E(C,B)$ is mapped to a $C$--$B$ edge
with probability $q/s^2$.  Hence the expected number of common
$C$--$B$ edges is
\[
 \frac{q^2}{s^2}.
\]
If $q<2s$, this expectation is less than four, so some choice gives at
most three common edges.

We may therefore assume $q=2s$.  Then every $b\in B$ has exactly two
neighbours in $C$.  Suppose that the first alternative never occurs.
The average above is four, so every pair of bijections gives exactly four
common edges.

Fix $\alpha$.  For $b\in B$ and $x\in C$, let
\[
 w_b(x)=\bigl|\{c\in N(b)\cap C:x\in N(\alpha(c))\cap C\}\bigr|.
\]
For every bijection $\beta$ the number of common $C$--$B$ edges is
\[
 \sum_{b\in B}w_b(\beta(b))=4.
\]
Fix distinct $b,b'\in B$ and distinct $x,y\in C$.  Choose a bijection
$\beta$ with $\beta(b)=x$ and $\beta(b')=y$, and let $\beta'$ be obtained
from $\beta$ by interchanging these two images.  Since both bijections
give exactly four common edges,
\[
 w_b(x)+w_{b'}(y)=w_b(y)+w_{b'}(x).
\]
Hence $w_b-w_{b'}$ is constant on $C$.  Since
\[
 \sum_{x\in C}w_b(x)=4
\]
for every $b$, this constant is zero, and all functions $w_b$ are equal.

Suppose that two vertices $b,b'$ have different two-element
neighbourhoods in $C$.  Then there is $x\in C$ with
$0<d(x,B)<s$.  Put
\[
 S_x(\alpha)=\{c\in C:x\in N(\alpha(c))\}
             =\alpha^{-1}(N(x)\cap B).
\]
As $\alpha$ varies, $S_x(\alpha)$ runs through all $d(x,B)$-subsets of
$C$, while
\[
 w_b(x)=|(N(b)\cap C)\cap S_x(\alpha)|.
\]
Thus $w_b=w_{b'}$ would imply that two distinct two-element subsets of
$C$ have the same intersection size with every $r$-subset, where
$0<r<s$, which is impossible.  Hence all vertices of $B$ have the same
two neighbours in $C$.
\cbdo

The exceptional configuration in Lemma~\ref{corebuffer} is broken by the
following degree identity.

\begin{lem}\label{defect}
Let $G$ be an $n$-vertex graph with $e(G)\leq2n-5$, and let
$P=\{p,q\}\subseteq V(G)$.  Then at least one of the following holds.
\begin{enumerate}
\item There is $v\notin P$ such that
\[
 d(v,P)\leq1\qquad\text{and}\qquad d(v)\leq2.
\]
\item There are at least
\[
 2+2\cdot\ind_{\{pq\in E(G)\}}
\]
vertices $v\notin P$ such that
\[
 d(v,P)=0\qquad\text{and}\qquad d_{G-P}(v)=3.
\]
\end{enumerate}
\end{lem}

\prf
We have
\[
 e(G)=2(n-2)-\sum_{v\notin P}(2-d(v,P))+e(G-P)
       +\ind_{\{pq\in E(G)\}}.
\]
The assumption $e(G)\leq2n-5$ therefore gives
\begin{equation}\label{defectsum}
 \sum_{v\notin P}
 \bigl(2(2-d(v,P))-d_{G-P}(v)\bigr)
 \geq2+2\cdot\ind_{\{pq\in E(G)\}}.
\end{equation}
If the first alternative fails, a vertex with $d(v,P)=1$ has
$d_{G-P}(v)\geq2$, while a vertex with $d(v,P)=0$ has
$d_{G-P}(v)\geq3$.  Hence the only positive summands in
\eqref{defectsum} come from vertices satisfying the second alternative,
and each contributes one. \cbdo

The following absorption argument is motivated by the three-cycle
construction used by Fang and Hou~\cite[Section~4.2]{FH}.  Here it is
adapted to a growing core and combined with our completion lemma, so the
degree estimates for the forbidden graph have to be established
separately.

\begin{lem}\label{absorb}
Let $G=G_n$ be an $n$-vertex graph with $e(G)\leq2n$, and let $C,B$ be
disjoint sets of the same size $c=c(n)$.  Put
\[
 D=\Delta(G-C).
\]
Suppose
\[
 c=o(n),\qquad D=o(n),\qquad cD=o(n),
\]
and suppose that
\begin{enumerate}
\item $B$ is independent, every $b\in B$ has degree at most three, and no
two vertices of $B$ have a common neighbour outside $C$;
\item $\alpha:C\to B$ and $\beta:B\to C$ are bijections which create at
most three common edges inside $C\cup B$;
\item there are constants $\gamma>0$ and $d_0$ such that, whenever
$N(b)\setminus C\neq\emptyset$, each of
$\alpha^{-1}(b)$ and $\beta(b)$ has at least $\gamma n$ non-neighbours
of degree at most $d_0$.
\end{enumerate}
Then, for all sufficiently large $n$, the partial permutation
\[
 c\mapsto\alpha(c)\quad(c\in C),
 \qquad
 b\mapsto\beta(b)\quad(b\in B)
\]
extends to a permutation of $V(G)$ with at most three common edges.
\end{lem}

\prf
Put
\[
 X=N(B)\setminus C.
\]
Then $|X|\leq3c$, and every $x\in X$ has a unique neighbour
$b(x)\in B$.  For each $x\in X$ choose two new vertices $p_x,q_x$, all
distinct and outside $C\cup B\cup X$, such that
\[
 d(p_x),d(q_x)\leq d_0,
\]
\[
 p_x\notin N(\beta(b(x))),
 \qquad
 q_x\notin N(\alpha^{-1}(b(x))),
\]
and, with
\[
 P_0=\{p_x:x\in X\},
 \qquad
 Q_0=\{q_x:x\in X\},
\]
the sets $P_0,Q_0$ are independent and
\[
 E(X,P_0)=E(X,Q_0)=E(P_0,Q_0)=\emptyset.
\]
The choice is greedy.  At every step start with at least $\gamma n$
low-degree non-neighbours of the required core vertex.  Excluding
$C\cup B\cup X$ removes $o(n)$ vertices, excluding $N_{G-C}(X)$ removes
at most $3cD=o(n)$ vertices, and previously chosen vertices together
with their neighbourhoods remove $O(cd_0)=o(n)$ vertices.  Thus the
choice is possible.

Extend the partial permutation by
\[
 x\mapsto p_x\mapsto q_x\mapsto x
 \qquad(x\in X).
\]
No new common edge is created among the prescribed vertices.  Edges
inside $C$ are sent into the independent set $B$; there are no edges
inside $B,P_0,Q_0$, and an edge inside $X$ is sent into $P_0$.  The
sets $X,P_0,Q_0$ are pairwise anticomplete as required above, and there
are no edges from $B$ to $P_0\cup Q_0$.  

It remains to check edges joining the core or the buffer to
$X\cup P_0\cup Q_0$.  An edge from $C$ to $X$ is mapped to a pair
between $B$ and $P_0$, and an edge from $C$ to $P_0$ to a pair between
$B$ and $Q_0$; neither can become common by the choice of
$P_0,Q_0$.  Next consider an edge $cq_x$ with $c\in C$.  Its image is
$\alpha(c)x$.  Since $x$ has the unique neighbour $b(x)$ in $B$, this
could be an edge only if $\alpha(c)=b(x)$, that is,
$c=\alpha^{-1}(b(x))$; this is excluded by the choice
$q_x\notin N(\alpha^{-1}(b(x)))$.

Finally, if $b(x)x$ is the edge joining $B$ to $X$, then its image is
$\beta(b(x))p_x$, which is a non-edge by the choice
$p_x\notin N(\beta(b(x)))$.  There are no edges from $B$ to
$P_0\cup Q_0$.  Hence no edge involving
$X\cup P_0\cup Q_0$ creates an additional common edge.
Thus only the at most three common edges already present inside
$C\cup B$ remain.

Let $p$ be the partial permutation defined on
\[
 S=C\cup B\cup X\cup P_0\cup Q_0.
\]
By the preceding verification, $I_G(p)\leq3$.  Since $p(S)=S$, the remaining
source and target set is
\[
 R=V(G)\setminus S.
\]
By Lemma~\ref{completion}, it remains only to list-pack the two copies of
$G[R]$ while avoiding the associated forbidden graph $Y_p$.

The remaining graph $G[R]$ has $O(n)$ edges and maximum degree at most
$D=o(n)$.  We claim that the forbidden graph $Y_p$ also satisfies
\[
 \Delta(Y_p)=o(n).
\]

Fix first a remaining source vertex $u\in R$.  Prescribed vertices in
$C\cup B$ create no forbidden targets: an image of a vertex of $C$
lies in $B$, which has no neighbours in $R$, while a vertex of $B$
itself has no neighbours in $R$.  The vertex $u$ has at most $D$
neighbours in each of $X$ and $P_0$, whose images lie in $P_0$ and
$Q_0$, respectively, and all vertices of $P_0\cup Q_0$ have degree at
most $d_0$.  These two classes therefore forbid at most
\[
 2d_0D
\]
targets.  Finally, $u$ has at most $|Q_0|=|X|\leq3c$ neighbours in
$Q_0$, and their images lie in $X$, where every vertex has at most
$D$ neighbours in $R$.  Thus
\[
 d_{Y_p}(u)\leq 2d_0D+3cD=o(n).
\]

The same estimate holds on the target side.  Indeed, prescribed pairs
involving $C\cup B$ again contribute nothing.  If $v\in R$ is adjacent
to a vertex of $P_0$, its preimage lies in $X$ and has at most $D$
neighbours in $R$; since $|P_0|\leq3c$, these pairs contribute at most
$3cD$.  If $v$ is adjacent to a vertex of $Q_0$ or $X$, the
corresponding preimage lies in $P_0$ or $Q_0$, respectively, and has
degree at most $d_0$.  Since $v$ has at most $D$ neighbours in each of
$Q_0$ and $X$, these two classes contribute at most $2d_0D$.  Hence
\[
 \Delta(Y_p)\leq 2d_0D+3cD=o(n),
\]
because $D=o(n)$ and $cD=o(n)$.

Since $|R|=n-o(n)$, Lemma~\ref{sparse} gives a list packing of the two
copies of $G[R]$ avoiding $Y_p$.  Lemma~\ref{completion} now extends
$p$ without increasing its number of common edges. \cbdo

We also record the simpler closed-buffer consequence.

\begin{lem}\label{closed}
Let $G=G_n$ be an $n$-vertex graph with $e(G)\leq2n$.  Let $C,B$ be
disjoint sets of the same size $c=o(n)$, and suppose that
\[
 E(B,V(G)\setminus(C\cup B))=\emptyset,
 \qquad
 \Delta(G-C)=o(n).
\]
Every permutation of $C\cup B$ which interchanges $C$ and $B$ and
creates at most three common edges extends, for all sufficiently large
$n$, to a permutation of $V(G)$ with at most three common edges.
\end{lem}

\prf
Let
\[
 R=V(G)\setminus(C\cup B)
\]
and let $p$ be the prescribed permutation of $S=C\cup B$.
The forbidden graph $Y_p$ from Lemma~\ref{completion} is empty.
Indeed, if $a\in C$ and $u\in R$ are adjacent, then
$p(a)\in B$, which has no neighbours in $R$.  On the other hand,
if $a\in B$, then $a$ itself has no neighbours in $R$.  Hence no
prescribed vertex can create a forbidden pair between the two
remaining copies of $R$.

Since $c=o(n)$,
\[
 |R|=n-2c=n-o(n).
\]
Moreover,
\[
 e(G[R])\leq2n=(2+o(1))|R|<3|R|
\]
for all sufficiently large $n$, and
\[
 \Delta(G[R])\leq\Delta(G-C)=o(n)=o(|R|).
\]
Thus Lemma~\ref{sparse} gives a list packing of the two copies of
$G[R]$ avoiding $Y_p=\emptyset$, which by Lemma~\ref{completion}
extends $p$ without creating any additional common edge. \cbdo

\section{Proof of Theorem~\ref{main}}

Suppose, for a contradiction, that there is an infinite sequence of
counterexamples $G=G_n$, of orders tending to infinity, such that
\begin{equation}\label{counter}
 e(G)\leq2n-5
 \qquad\text{and}\qquad
 I_G(\sigma)\geq4\quad\text{for every }\sigma.
\end{equation}
By Lemma~\ref{sparse}, applied with no forbidden pairs, there is a fixed
$\varepsilon>0$ such that
\begin{equation}\label{linearDelta}
 \Delta(G)\geq\varepsilon n
\end{equation}
for all sufficiently large graphs in the sequence.

After passing to a subsequence, one of the following holds:
\begin{enumerate}
\item there is a constant $\gamma>0$ such that
$\Delta(G)\leq(1-\gamma)n$;
\item a vertex of maximum degree has $o(n)$ non-neighbours.
\end{enumerate}
We treat these cases separately.

\subsection{No almost-universal vertex}

The organization of this case follows the core--buffer strategy of
Fang and Hou~\cite[Section~4]{FH}, with two changes that are essential
here: our core is allowed to grow, and the equality case in the
first-moment bound must be resolved rather than accepted as four common
edges.

Assume
\begin{equation}\label{nonuniversal}
 \Delta(G)\leq(1-\gamma)n
\end{equation}
for some fixed $\gamma>0$.  Apply Lemma~\ref{corechoice} with
$m=\lfloor n^{1/3}\rfloor$, and retain the notation $C,c,D$ from that
lemma.  By \eqref{linearDelta},
\begin{equation}\label{linearCoreDegree}
 S_C:=\sum_{v\in C}d(v)\geq\varepsilon n.
\end{equation}
Lemma~\ref{reservoir} therefore gives
\begin{equation}\label{linearBcal}
 |\mathcal B|=\Omega(n),
\end{equation}
where
\[
 \mathcal B=\{v\notin C:d(v)\leq3,\,d(v,C)\leq2\}.
\]

Choose $B\subseteq\mathcal B$, $|B|=c$, so that $B$ is independent and
no two vertices of $B$ have a common neighbour outside $C$.  This is the
greedy buffer selection from Fang and Hou~\cite[Section~4.1]{FH}; the
same local count applies here.  Once $b$ is chosen, at most $3D+4$
candidates are lost, while
\[
 c(3D+4)=o(n)
\]
by Lemma~\ref{corechoice} and $\mathcal B$ is linear.

Every vertex has at least $\gamma n-1$ non-neighbours.  Choose a fixed
$d_0$ so large that $4/d_0<\gamma/2$.  At most $4n/d_0$ vertices have
degree greater than $d_0$, and hence every vertex has at least
$\gamma n/3$ non-neighbours of degree at most $d_0$ for large $n$.

Apply Lemma~\ref{corebuffer}.  If its first alternative holds, the
resulting bijections $\alpha,\beta$ satisfy the assumptions of
Lemma~\ref{absorb}, which gives a permutation with at most three common
edges, contradicting \eqref{counter}.

It remains to consider the rigid alternative.  We first show that there
is a fixed pair $P\subseteq C$ such that
\begin{equation}\label{allP}
 N(b)\cap C=P\qquad\text{for every }b\in\mathcal B.
\end{equation}

Call two vertices of $\mathcal B$ compatible if they are nonadjacent and
have no common neighbour outside $C$, and let $K$ be the graph on
$\mathcal B$ in which two vertices are adjacent precisely when they are
compatible.  Every vertex is incompatible with at most $3D+4=o(n)$
other vertices.  Since $|\mathcal B|=\Omega(n)$, for all sufficiently
large $n$ we have $3D+4<|\mathcal B|/2$, and hence
\[
 \delta(K)>|\mathcal B|/2-1.
\]
In particular, $K$ is connected.  Moreover, every compatible pair
extends greedily to a compatible $c$-set, because $cD=o(n)$ and
$|\mathcal B|=\Omega(n)$.  If a compatible pair had different
neighbourhoods in $C$, extend it to such a set $B$.  The second
alternative of Lemma~\ref{corebuffer} would then be impossible, while
the first, followed by Lemma~\ref{absorb}, would contradict
\eqref{counter}.  Hence adjacent vertices of $K$ have the same
neighbourhood in $C$.  Since $K$ is connected, for any
$b,b'\in\mathcal B$ there is a path
$b=b_0,b_1,\ldots,b_\ell=b'$ in $K$; equality of the $C$-neighbourhoods
along this path gives
$N(b)\cap C=N(b')\cap C$.  This proves \eqref{allP}.

Apply Lemma~\ref{defect} to $P=\{p,q\}$.  Its first alternative is
impossible.  Indeed, $C$ contains only vertices of degree at least four,
while a vertex $v\notin C$ with $d(v)\leq2$ and $d(v,P)\leq1$ belongs
to $\mathcal B$ but does not satisfy \eqref{allP}.

Consequently there are at least two vertices $y,z\notin C$ such that
\[
 d(y,P)=d(z,P)=0,
 \qquad
 d_{G-P}(y)=d_{G-P}(z)=3.
\]
They have total degree three.  Since neither can belong to $\mathcal B$,
all three neighbours of each lie in $C$.  In particular $y$ and $z$ are
nonadjacent.  Put
\[
 A=N(y),\qquad D_0=N(z).
\]
Thus $A,D_0\subseteq C\setminus P$ and $|A|=|D_0|=3$.
Choose a new compatible buffer consisting of $y,z$ and $c-2$ vertices
of $\mathcal B$.  The latter all have neighbourhood $P$ in $C$.

If $A\cap D_0\neq\emptyset$, choose distinct $r,s\in A\cup D_0$ such
that
\[
 \ind_{\{r\in A\}}+\ind_{\{r\in D_0\}}+
 \ind_{\{s\in A\}}+\ind_{\{s\in D_0\}}\geq3.
\]
Set
\[
 \alpha(r)=y,
 \qquad
 \alpha(s)=z,
 \qquad
 \beta(y)=p,
 \qquad
 \beta(z)=q,
\]
and complete both bijections using the vertices with neighbourhood $P$.
The only common $C$--$B$ edges are among the edges from
$A\setminus\{r,s\}$ to $y$ and from $D_0\setminus\{r,s\}$ to $z$.
There are at most three.

If $A\cap D_0=\emptyset$, write $D_0=\{d_1,d_2,d_3\}$.  Set
\[
 \alpha(d_2)=y,
 \qquad
 \alpha(d_3)=z,
 \qquad
 \beta(y)=d_1,
 \qquad
 \beta(z)=q,
\]
and send one of the remaining buffer vertices to $p$ under $\beta$.
Complete both bijections arbitrarily.  One common edge is incident with
$z$, and two are incident with the buffer vertex sent to $p$; no other
common edge occurs.  Thus again there are at most three.

In either case the two exceptional vertices have no neighbour outside
$C$, while all other buffer vertices satisfy the hypotheses of
Lemma~\ref{absorb}.  This gives the final contradiction in the case
\eqref{nonuniversal}.

\subsection{An almost-universal vertex}

We may now assume that a vertex $c$ of maximum degree has
\[
 d(c)=n-1-t,
 \qquad
 t=o(n).
\]
Let $d_2$ be the second largest degree.  After passing to a further
subsequence, either $d_2=n-o(n)$ or there is a constant $\gamma>0$ such
that
\begin{equation}\label{secondgap}
 d_2\leq(1-\gamma)n.
\end{equation}

\subsubsection{Two almost-universal vertices}

Assume first that two vertices $p,q$ have degree $n-o(n)$, and put
$P=\{p,q\}$.  Let
\[
 M=\sum_{v\notin P}(2-d(v,P)).
\]
Then $M=o(n)$.  Since
\[
 e(G)=2(n-2)-M+e(G-P)+\ind_{\{pq\in E(G)\}},
\]
the inequality $e(G)\leq2n-5$ gives $e(G-P)=o(n)$.  Hence all but
$o(n)$ vertices outside $P$ satisfy
\begin{equation}\label{trueTwins}
 N(v)=P.
\end{equation}
We call them $P$-twins.  Notice also that
$\Delta(G-P)=o(n)$.

Apply Lemma~\ref{defect}.  Suppose first that there is $v\notin P$ with
$d(v,P)\leq1$ and $d(v)\leq2$.  Put
\[
 C=P\cup(N(v)\setminus P)
\]
and let $B$ consist of $v$ and $|C|-1$ $P$-twins.  Then
$E(B,V(G)\setminus(C\cup B))=\emptyset$, and the neighbourhoods in $C$
have size at most two and are not all equal to $P$.  Lemma~\ref{corebuffer}
and then Lemma~\ref{closed} give a contradiction.

We may therefore assume the second alternative of Lemma~\ref{defect}.
If two of its vertices $y,z$ are nonadjacent, put
\[
 C=P\cup N(y)\cup N(z)
\]
and fill $B$ with $y,z$ and $P$-twins.  With
$A=N(y)$ and $D_0=N(z)$, the two constructions used in the preceding
subsection give bijections with at most three common edges.  The buffer
is closed, so Lemma~\ref{closed} completes the permutation.

It remains that all vertices supplied by the second alternative are
pairwise adjacent.  Choose two of them, $y,z$, and put
\[
 C=P\cup(N(y)\setminus\{z\})\cup(N(z)\setminus\{y\}).
\]
Write
\[
 A=N(y)\cap C,
 \qquad
 D_0=N(z)\cap C.
\]
Then $|A|=|D_0|=2$.  Fill $B$ with $y,z$ and $P$-twins.  Choose
distinct $r\in A$ and $s\in D_0$, and set
\[
 \alpha(r)=y,
 \qquad
 \alpha(s)=z,
 \qquad
 \beta(y)=p,
 \qquad
 \beta(z)=q,
\]
completing both bijections by the $P$-twins.  The number of common edges
inside $C\cup B$ is at most
\begin{equation}\label{connectedcount}
 4-
 \bigl(\ind_{\{r\in A\}}+\ind_{\{r\in D_0\}}+
       \ind_{\{s\in A\}}+\ind_{\{s\in D_0\}}\bigr)
 +\ind_{\{rs\in E(G)\}}+
 \ind_{\{pq\in E(G)\}}.
\end{equation}
If $pq\notin E(G)$, this is at most three.  If $pq\in E(G)$, the defect
lemma supplies at least four such vertices.  Since they are pairwise
adjacent and each has degree three in $G-P$, they form a $K_4$.  For
any adjacent pair $y,z$, the other two vertices form both $A$ and $D_0$;
choosing them for $r,s$ makes \eqref{connectedcount} equal to two.
The generalized closed-buffer lemma applies even though $yz$ is an edge,
and gives a contradiction.

\subsubsection{Exactly one almost-universal vertex}

This part refines the almost-universal argument of Fang and
Hou~\cite[Lemma~3.1]{FH}.  Their construction permits four common
edges; the prescribed images below are chosen so that one further
possible common edge is eliminated.

Assume \eqref{secondgap}.  Put $J=G-c$ and
\[
 L=\{v\in V(J):d_J(v)\leq1\}.
\]
We first consider $|L|=o(n)$.  Since
\[
 e(J)\leq n-4+t,
\]
we have
\[
 \sum_{v\in V(J)}(d_J(v)-2)
 =2e(J)-2(n-1)
 \leq 2t-6.
\]
Vertices of degree two contribute nothing to this sum, while vertices
in $L$ contribute negatively.  Hence
\[
 \sum_{\substack{v\in V(J)\\ d_J(v)\geq3}}(d_J(v)-2)
 \leq
 2t-6+\sum_{v\in L}(2-d_J(v))
 \leq 2t-6+2|L|.
\]
Therefore
\begin{equation}\label{excessJ}
 \sum_{\substack{v\in V(J)\\ d_J(v)\geq3}}(d_J(v)-2)
 \leq2t-6+2|L|=o(n).
\end{equation}
Consequently, with
\[
 H=\{v\in V(J):d_J(v)>2\},
\]
we have
\begin{equation}\label{HJ}
 |H|=o(n),
 \qquad
 \sum_{v\in H}d_J(v)=o(n),
 \qquad
 \Delta(J)=o(n).
\end{equation}
Since the left-hand side of \eqref{excessJ} is nonnegative,
\[
 |L|\geq 3-t.
\]
Thus the choices of a vertex from $L$ made below for $t=0,1$ are possible.

If $t=0$, choose $a\in L$ and put $F=G-\{c,a\}$.  Then
$e(F)\leq|V(F)|-2$, so $F$ has a proper self-packing.  Extend it by
swapping $c$ and $a$.  The edge $ca$ contributes at most one common
edge, and the at most one edge from $a$ into $F$ can contribute once in
each direction.  Thus there are at most three common edges.

If $t=1$, again choose $a\in L$.  If $d_J(a)=1$, the preceding proper
packing argument applies.  If $d_J(a)=0$, then
$e(F)\leq|V(F)|-1$; by Proposition~\ref{firstmoment} with $k=1$, $F$ has a relabelling
with at most one common edge.  Swapping $c,a$ adds at most one further
common edge.

Let $t\geq2$.  By \eqref{excessJ}, all but $o(n)$ vertices of $J$ have
degree two.  Let $T=V(J)\setminus N_G(c)$.  Choose a degree-two vertex
\[
 a\notin T\cup H\cup N_J(H)\cup N_J(T\setminus H).
\]
This is possible because all excluded sets have size $o(n)$.  Write
\[
 N_J(a)=\{z_1,z_2\}.
\]
Then $z_1,z_2\notin T\cup H$.  Since $z_1$ has degree at most two and
is adjacent to $a$, it has at most one neighbour in $T$.  Choose
distinct $b_1,b_2\in T$ such that
\[
 b_2z_1\notin E(G).
\]
On $F=G-\{c,a\}$ prescribe
\[
 z_2\mapsto b_2,
 \qquad
 b_1\mapsto z_1.
\]
This partial bijection has no common edge, because the only possible source
edge between the prescribed source vertices is sent to the non-edge
$b_2z_1$.  Lemma~\ref{prescribed}, using \eqref{HJ}, extends it to a
proper self-packing $\phi$ of $F$.

Now swap $c$ and $a$ and use $\phi$ on $F$.  The edge $ca$ contributes
one common edge.  Of the two edges $az_1,az_2$, the second is sent to
$cb_2$, a non-edge, so at most one survives.  In the reverse direction,
a $c$-edge can become an $a$-edge only through a preimage of $z_1$ or
$z_2$; the preimage of $z_1$ is $b_1$, and $cb_1$ is a non-edge.  Thus
this direction contributes at most one further common edge.  The total
is at most three.

We are left with $|L|=\Omega(n)$.  Since $|T|=o(n)$, all but $o(n)$
vertices of $L$ are adjacent to $c$.  After passing to a subsequence,
one of the sets
\[
 L_0=\{v:d_J(v)=0,\ cv\in E(G)\},
 \qquad
 L_1=\{v:d_J(v)=1,\ cv\in E(G)\}
\]
is linear.

Suppose first that $|L_0|=\Omega(n)$.  Apply Lemma~\ref{corechoice} with
$m=\lfloor n^{1/3}\rfloor$ and let $C$ be the resulting core; it
contains $c$.  
Take $|C|$ vertices of $L_0$ for $B$.  Since every $b\in B$ satisfies
$N_G(b)=\{c\}\subseteq C$, we have
\[
 E(B,V(G)\setminus C)=\emptyset.
\]
Under bijections $C\to B$ and $B\to C$ there is exactly one common
$C$--$B$ edge: it is
the edge whose buffer endpoint is mapped to $c$.  Lemma~\ref{closed}
gives a contradiction.

It remains that $|L_1|=\Omega(n)$.  For $b\in L_1$ let $z(b)$ be its
unique neighbour in $J$, so
\begin{equation}\label{leafpair}
 N_G(b)=\{c,z(b)\}.
\end{equation}
Pass to a subsequence on which either some value $z$ occurs for
$\Omega(n)$ vertices $b$, or the maximum multiplicity of a value $z(b)$
is $o(n)$.

In the first case put $P=\{c,z\}$.  Then linearly many vertices
$b\in L_1$ satisfy $N_G(b)=P$, and hence are $P$-twins.
Apply Lemma~\ref{corechoice} with
$m=\lfloor n^{1/3}\rfloor$ and enlarge the resulting core, if necessary,
to contain $P$ and the bounded exceptional neighbourhoods used below.
The enlarged core still has order $o(n)$ and leaves maximum degree
$o(n)$ outside it.

If the first alternative of Lemma~\ref{defect} gives a vertex $v$, add
$N(v)\setminus P$ to the core and take for $B$ the vertex $v$ together
with $P$-twins.  Lemmas~\ref{corebuffer} and~\ref{closed} give a
contradiction.  Otherwise use vertices from the second alternative.  If
two are nonadjacent, add their neighbours to the core and use the two
three-neighbour constructions from the non-almost-universal case.  If
all are pairwise adjacent, use the construction leading to
\eqref{connectedcount}. 
In every case
\[
 E(B,V(G)\setminus(C\cup B))=\emptyset,
\]
and the prescribed permutation of $C\cup B$ creates at most three
common edges.  Hence Lemma~\ref{closed} completes the permutation.

We may therefore assume that the maximum multiplicity of $z(b)$ is
$o(n)$.  Let $r$ be the number of distinct values $z(b)$.  Then
$r\to\infty$.  Apply Lemma~\ref{corechoice} with
\[
 m=\min\left\{\left\lfloor\frac r4\right\rfloor,
               \left\lfloor n^{1/3}\right\rfloor\right\}.
\]
Let $C$ be the resulting core.  Choose $y\in L_1$ with $z(y)\notin C$
and add
\[
 z=z(y)
\]
to $C$.  The enlarged core still satisfies
\[
 |C|=o(n),
 \qquad
 |C|\Delta(G-C)=o(n).
\]
Choose $B$, $|B|=|C|$, consisting of $y$ and further vertices of $L_1$,
so that the values $z(b)$ are distinct, all $z(b)$ with $b\neq y$ lie
outside $C$, and $B$ is independent.  This is possible because the
number of available distinct values is more than twice $|C|$, while a
chosen leaf is adjacent to at most one other leaf.

Choose $b_0\in B\setminus\{y\}$ and define
\[
 \alpha(c)=y,
 \qquad
 \alpha(z)=b_0,
 \qquad
 \beta(y)=c,
 \qquad
 \beta(b_0)=z,
\]
with
\[
 \alpha(C\setminus\{c,z\})=B\setminus\{y,b_0\},
 \qquad
 \beta(B\setminus\{y,b_0\})=C\setminus\{c,z\}.
\]
There are exactly three common edges inside $C\cup B$:
\[
 cy,\qquad zy,\qquad cb_0.
\]
The vertex $y$ has no neighbour outside $C$.  Every other buffer vertex
which has an outside neighbour has both its $\alpha$-preimage and its
$\beta$-image different from $c$.  By \eqref{secondgap}, every such core
vertex has at least $\gamma n/3$ non-neighbours of bounded degree, after
choosing the degree bound sufficiently large.  Lemma~\ref{absorb}
therefore completes the permutation, a contradiction.

All possible subsequences of counterexamples have been excluded.  Hence,
for all sufficiently large $n$, every $n$-vertex graph with at most
$2n-5$ edges has an $\mathcal E_3$-near-packing.  By
\eqref{relation-f-mu} and the construction $K_{2,n-2}$,
\[
 f(n,4)=2n-4.
\]
\cbdo

\section*{Declaration on the use of artificial intelligence}
Generative artificial intelligence tools were used during the
preparation of this manuscript to improve the exposition and to assist
in developing and drafting parts of the proof.  In particular, several
auxiliary arguments in the proof of Theorem~\ref{main} were initially
developed with the assistance of GPT-5.6.  All mathematical statements
and proofs presented here have been checked and revised by the author,
who assumes full responsibility for their correctness.

\end{document}